\documentclass[11pt]{article}\include{ddefs}
\usepackage{graphics,pifont,array}
\usepackage{graphicx}
\usepackage{endnotes}
\let\footnote=\endnote

\begin{document}

\begin{titlepage}

\begin{center}

{\Large{        \mbox{   }                          \\
                \mbox{   }                          \\
                \mbox{   }                          \\
                \mbox{   }                          \\
                \mbox{   }                          \\
                \mbox{   }                          \\
                \mbox{   }                          \\
                \mbox{   }                          \\
                \mbox{   }                          \\
                \mbox{   }                          \\
                \mbox{   }                          \\
  {\textbf{How {\em{Not}} to Compute a Fourier Transform}}  \\
               \mbox{    }                          \\
               \mbox{    }                          \\
              J. A. Grzesik                         \\
           Allwave Corporation                      \\
        3860 Del Amo Boulevard                      \\
                Suite 404                           \\
           Torrance, CA 90503                       \\
                \mbox{    }                         \\
        +1 (818) 749-3602                           \\ 
            jan.grzesik@hotmail.com                 \\
              \mbox{     }                          \\
              \mbox{     }                          \\  
              \today                                      }  }

\end{center}

\end{titlepage}

\pagenumbering{roman}

\setcounter{page}{2}

\vspace*{+2.725in}

\begin{abstract}

	\parindent=0.245in
	
	We revisit the Fourier transform of a Hankel function, of considerable importance
	in the theory of knife edge	diffraction.  Our approach is based directly upon the
	underlying Bessel equation, which admits manipulation into an alternate second
	order differential equation, one of whose solutions is precisely the desired
	transform, apart from an {\em{a priori}} unknown constant and a second, undesired
	solution of logarithmic type.  A modest amount of analysis is then required to
	exhibit that constant as having its proper value, and to purge the logarithmic
	accompaniment.  The intervention of this analysis, which relies upon an interplay
	of asymptotic and close-in functional behaviors, prompts our somewhat ironic,
	mildly puckish caveat,
	our negation {\em{``not''}} in the title.  In a concluding section we show
	that this same transform is still more readily exhibited as an easy byproduct of
	the inhomogeneous wave equation in two dimensions satisfied by the Green's
	function $G,$ itself proportional to a Hankel function.  This latter discussion
	lapses of course into the argot of physicists and, in its r{\^{o}}le of a mere
	afterthought, makes no claim whatsoever to any originality.
	\parindent=0.0in
	\newline

	{\bf{key words}} --- knife edge diffraction; Fourier cosine transform; Bessel's
	differential equation; asymptotic Hankel function estimate; two-dimensional Laplace equation;
	Green's	function for two-dimensional Helmholtz equation

\end{abstract}

\newpage

\parindent=0.5in

\pagenumbering{arabic}

\setcounter{page}{1}

\pagestyle{myheadings}

\mbox{   }

\markright{J. A. Grzesik \\ how not to compute a fourier transform}

\section{Introduction}
\vspace{-4mm}

   Sommerfeld's celebrated, period $4\pi$ treatment of knife edge diffraction [{\bf{1}},{\bf{2}}] was quick to elicit
alternate viewpoints [{\bf{3}},{\bf{4}}] of comparable ingenuity if, regrettably, far less renown.  These achievements
rested on their deserved laurels for half a century, only to be roused from their slumbers during the early forties
by a minor renaissance anchored around the Wiener-Hopf [{\bf{W-H}}] technique, first reported in [{\bf{5}},{\bf{6}}],
and then belatedly in [{\bf{7}}].  Although initiated in a framework of integral equations, the {\bf{W-H}} program soon
abandoned this mooring in favor of Fourier transforms applied directly to the underlying differential (wave) equation
[{\bf{8}},{\bf{9}},{\bf{10}}], and this latter approach, too, was ultimately complemented by the method of dual
integral equations {\em{cum}} plane wave spectrum [{\bf{11}},{\bf{12}},{\bf{13}}].

    The integral equation springboard in [{\bf{5}},{\bf{6}},{\bf{7}}] pivots around the Green's function for the
two-dimensional wave equation, a quantity proportional to the Hankel function $H_{0}^{(2)}(k\rho)$ of the
second kind\footnote{A Hankel function of the second kind is paired with the adoption of a simple harmonic dependence
upon time $t$ evolving in accordance with $\exp(i\omega t),\;\omega = 2\pi/T > 0,$ with $T$ being the period.  It is likewise
assumed that $\Re k > 0,\;\Im k \leq 0,$ the latter intended to account for the physically welcome possibility of wave
dissipation.  In the absence thereof, $k=2\pi/\lambda,$ $\lambda = cT$ being the wavelength dictated by light speed $c.$
Symbol $\rho = \sqrt{x^{2}+y^{2}\,}$ measures the two-dimensional distance from co\"{o}rdinate origin.} and order zero.
There quickly ensues a need for the Fourier transform thereof having the known outcome\footnote{While our notation aims
to be as close as possible to that of [{\bf{5}},{\bf{6}}], we have dispensed with an extraneous factor of $\sqrt{2/\pi\,},$
and, with a view to its imminent, analytically more robust r\^{o}le, have included wave number $k$ as an explicit argument.
Functional form $1/\sqrt{k^{2}-w^{2}\,}$ in line two, crucial for the succeeding {\bf{W-H}} developments, is gotten under
analytic continuation from a formula due to Alfred Barnard Basset [{\bf{14}}, p. 32], but stated only as an exercise no less,
without any indication of an overt proof (!!), and subsequently quoted in [{\bf{15}}, p. 388, Eq. (10)].  Details of the
requisite analytic continuation are provided in [{\bf{5}}, pp. 25-26].  The branch of $\sqrt{k^{2}-w^{2}\,}$ is understood
to give $k$ when $w=0.$
\newpage
\mbox{   }
\newline}   
\begin{eqnarray}
L(k,w) & = & \int_{0}^{\,\infty}H_{0}^{(2)}(kx)\cos(wx)dx \nonumber  \\
     & = & \frac{1}{\sqrt{k^{2}-w^{2}\,}}\,.
\end{eqnarray}
\vspace{-7mm}    

    Our aim here is to arrive at (1) by an alternate route, one which exploits the underlying Bessel
equation satisfied by $H_{0}^{(2)}(kx),$ with an emphasis upon wave number $k$ as variable and with $x$
temporarily regarded as an incidental multiplier.  We note of course that the quadrature in (1) is valid
throughout the horizontal strip $|\Im w | < |\Im k|,$ which allows us to restrict initial attention to       	
the situation with $\Im w = 0,\;\Re w > 0,$ and then to extend its outcome throughout that
strip\footnote{Indeed, such extension can sweep beyond this strip so as to include the entire $w$
plane, suitably slit with branch cuts emanating from $w_{\pm}=\pm k.$}
on the basis of analytic continuation.  And then, since incidental multiplier $x$ in Bessel's equation is
readily produced via differentiation with respect to $w,$ the ordinary differential equation (ODE) is
promoted at once to partial differential equation (PDE) status in variable pair $(k,w).$

    On exploiting next a variable scaling $\zeta=k/w,$
\begin{eqnarray}
L(k,w) & = & \frac{1}{w}\int_{0}^{\,\infty}H_{0}^{(2)}(kx/w)\cos(x)dx \nonumber  \\
& = & \frac{1}{w}F(\zeta)\,,
\end{eqnarray}
our PDE regresses to second-order ODE status in $\zeta,$ an equation among whose two independent, easily
ascertained solutions one indeed finds the second line from (1), but with an initially undetermined constant
multiplier, and accompanied by a somewhat more complicated, unwelcome logarithmic structure, similarly prefaced
with its own constant factor.

     And so, after this initial burst of easy progress, the remainder of the calculation devolves into
fixing these constants in a way to dismiss that logarithmic companion and indeed to match
(1).  This latter activity, which marries the asymptotic behavior of $L$ to the close-in behavior of
$H^{(2)}$ {\em{vis-\`{a}-vis}} their respective arguments, requires a modicum of spade work, diligence
which underwrites our tongue-in-cheek qualifier {\em{``how not to compute''}} embedded in note title.
On the other hand, diligent or not, this \mbox{latter effort,}
\newpage
\mbox{   }
\newline
\newline        
\newline
too, is not without some minor flecks of
mathematical allure.  Such, at any rate, is our own modest perception thereof.

    We conclude the note with yet another derivation of (1), couched in the physicist's language of
Green's functions and Dirac deltas.  Hopefully the tone of this grace note will not ruffle too many
mathematical sensibilities.  This material is essentially common folklore among physicists and is included,
bereft of any claim to originality, mainly for the sake of completeness.  We beg editorial indulgence
for allowing us thus to gild the lily.
\vspace{-5mm}
\section{From Bessel's ODE to PDE and then back to an ODE once more}
\vspace{-2mm}
   Bessel's ODE equation for $H_{0}^{(2)}(kx)$ yields in the first instance
\begin{eqnarray}	
\left[\rule{0mm}{7mm}k\frac{\partial^{2}}{\partial k^{2}}+
\frac{\partial}{\partial k}+x^{2}k\right]H_{0}^{(2)}(kx) & = & 0\,,
\end{eqnarray}               
whence it readily follows that          
\begin{eqnarray}
\left[\rule{0mm}{7mm}k\frac{\partial^{2}}{\partial k^{2}}+
\frac{\partial}{\partial k}-k\frac{\partial^{2}}{\partial w^{2}}\right]L(k,w) & = & 0\,,
\end{eqnarray}
Passage from (3) to (4) is premised on the assumption that the differential operator
\[k\frac{\partial^{2}}{\partial k^{2}} + \frac{\partial}{\partial k} - k\frac{\partial^{2}}{\partial w^{2}}\]
can be transferred with analytic impunity, {\emph{en masse}}, to the interior of the defining quadrature
from (1). This should certainly be legitimate since, even though its terms destroy convergence when
permitted to operate on $H_{0}^{(2)}(kx)\cos (wx)$ individually, in unison they simply annihilate
$H_{0}^{(2)}(kx)\cos (wx),$ yielding a constant, null outcome about whose integrability there can assuredly be
no doubt.
  
   Introduction of (2) has the effect of converting
PDE (4) into the homogeneous ODE\footnote{A subscript, such as $\zeta,$ denotes a derivative with respect to
the variable indicated.}
\begin{eqnarray}
(\zeta-\zeta^{3})F_{\zeta\zeta}(\zeta)+(1-4\zeta^{2})F_{\zeta}(\zeta)-2\zeta F(\zeta) & = & 0
\end{eqnarray}
which further condenses into
\begin{eqnarray}
\left(\rule{0mm}{4mm}(\zeta-\zeta^{3})F(\zeta)\right)_{\zeta\zeta}+
\left(\rule{0mm}{4mm}(2\zeta^{2}-1)F(\zeta)\right)_{\zeta} & = & 0
\end{eqnarray}
and hence allows an immediate integration into the form
\begin{eqnarray}
F(\zeta) & = & \frac{A}{\sqrt{1-\zeta^{2}\,}}+
\frac{B}{\sqrt{1-\zeta^{2}\,}}\log\left(\frac{\zeta}{1+\sqrt{1-\zeta^{2}\,}}\right)\,,
\end{eqnarray}
with constants $A$ and $B$ still to be determined.  And, from what has been said before, the desired form (1)
now literally drops into our lap, provided only that values $A=i$ and $B=0$ can somehow be assured.
\newpage
\mbox{   }
\section{Fixing integration constants {\mbox{\boldmath$A$}} and {\mbox{\boldmath$B$}}}
Constants $A$ and $B$ are found next by utilizing the dominant terms of $H_{0}^{(2)}(kx/w)$ as $w\rightarrow\infty.$
Thus, from (2) and (7) it follows on the one hand that 
\begin{eqnarray}
L(k,w) & \approx{\atop_{\rule{-6.5mm}{2.5mm}w\rightarrow\infty}} &  \frac{A}{w}+\frac{B}{w}\log\left(\frac{k}{2w}\right)\,. 
\end{eqnarray}
At the same time, from the first line of (2),
\begin{eqnarray}
L(k,w) & = \;\; &  \lim{\atop_{\rule{-7.5mm}{2.5mm}\beta\,\downarrow\,0+}}
        \frac{1}{w}\int_{0}^{\,\infty}e^{-\beta x}\cos(x)H_{0}^{(2)}(kx/w)dx   \nonumber   \\
       & \approx{\atop_{\rule{-6.5mm}{2.5mm}w\rightarrow\infty}} &
        \lim{\atop_{\rule{-7.5mm}{2.5mm}\beta\,\downarrow\,0+}}
       \frac{1}{w}\int_{0}^{\,\infty}e^{-\beta x}\cos(x)\left[\rule{0mm}{7mm}1-
           \frac{2i}{\pi}\left\{\rule{0mm}{6mm}\log\left(\frac{k}{2w}\right)+\gamma\right\}-
           \frac{2i}{\pi}\log(x)\right]dx \,,
\end{eqnarray}
with $\gamma=0.5772156649\ldots$ being the Euler-Mascheroni constant.\footnote{Our abrupt,
{\em{deus ex machina}} recourse, in the second line of (9), to the lowest order only of the
close-in Hankel function behavior is triggered by the circumstance that all of its remaining
terms are burdened by still higher powers of $1/w$ [{\bf{16}}].  We indicate in the Remarks
section at note's end how some of these terms can be disposed of.}  Here
\begin{eqnarray}
\int_{0}^{\infty}e^{-\beta x}\cos(x)dx & = & \frac{\beta}{1+\beta^{2}}
\end{eqnarray}
and thus vanishes together with $\beta,$ something which allows us to say at once that $B=0.$

     It remains hence to contend with
\begin{eqnarray}
\int_{0}^{\,\infty}e^{-\beta x}\cos(x)\log(x)dx & = & \beta\int_{0}^{\,\infty}e^{-\beta x }\sin(x)\log(x)dx
                        -\int_{0}^{\,\infty}e^{-\beta x} \frac{\sin(x)}{x}dx\,,
\end{eqnarray}      
its right-hand side gotten under a routine integration by parts.  And, while the second integral on the right
has of course the known limit $\pi/2$ when $\beta\downarrow 0+,$ a short excursus intrudes now in order to
ascertain the corresponding limit for the first.

    We begin by writing
\begin{eqnarray}
\beta\int_{0}^{\,\infty}e^{-\beta x }\sin(x)\log(x)dx & = & \beta M(\beta,1)\,,
\end{eqnarray}
with function $M(\beta,\xi)$ defined as
\begin{eqnarray}
M(\beta,\xi) & = & \int_{0}^{\,\infty}e^{-\beta x }\sin(\xi x)\log(x)dx \,.
\end{eqnarray}
It is evident by inspection that $M(\beta,\xi)$ satisfies the Laplace equation
\begin{eqnarray}
M(\beta,\xi)_{\beta\,\beta}+M(\beta,\xi)_{\xi\,\xi} & = & 0 \,.
\end{eqnarray}
On the other hand, at least in the first quadrant with $\beta > 0,$ $\xi > 0,$
\begin{eqnarray}
M(\beta,\xi) & = & \xi^{-1}M(\beta/\xi,1)-\xi^{-1}\log(\xi)\int_{0}^{\,\infty}e^{-\beta x/\xi}\sin(x)dx \nonumber  \\
    & = & \xi^{-1}M(\beta/\xi,1)-\xi\log(\xi)/(\beta^{2}+\xi^{2})\,.
\end{eqnarray}
\newpage
\mbox{   }
\newline
\newline        
\newline
By direct calculation\footnote{In both (16) and (17) hand manipulations have been successfully collated against their
{\bf{MATHEMATICA}} counterparts.}	
Laplace's equation (14) is then reduced to the following inhomogeneous ODE
\begin{eqnarray}
\left(\rule{0mm}{4mm}(1+\eta^{2})N(\eta)\right)_{\eta\,\eta} & = & (1+\eta^{2})^{-1}-4(1+\eta^{2})^{-2}
\end{eqnarray}
for the function $N(\eta)=M(\eta,1).$  Integration gives
\begin{eqnarray}
N(\eta) & = & \frac{1}{1+\eta^{2}}\left\{\rule{0mm}{6mm}C+D\eta-\frac{1}{2}\log(1+\eta^{2})-\eta\arctan(\eta)\right\}
\end{eqnarray}
with still other constants $C$ and $D$ whose precise values are not of any immediate interest.\footnote{Additional comments
regarding $C$ and $D,$ but short of any immediate evaluation, are likewise found in a concluding Remarks section.}
What is of interest is our newly acquired ability to state with some confidence that
\begin{equation}
\lim{\atop_{\rule{-7.5mm}{2.5mm}\beta\,\downarrow\,0+}}\beta M(\beta,1)=
\lim{\atop_{\rule{-7.5mm}{2.5mm}\beta\,\downarrow\,0+}}\beta N(\beta)=0\,.
\end{equation}

     And finally, on putting together the remnant of (9) with that of (11), we duly find
\begin{eqnarray}
L(k,w) & \approx{\atop_{\rule{-6.5mm}{2.5mm}w\rightarrow\infty}} & \frac{i}{w} \,,
\end{eqnarray}
so that indeed $A=i$ and we are done.
\section{Green's function byproduct}
We would be highly remiss were we not to note that the simplest, if perhaps not the most mathematically palatable
access route to transform (1) originates with the PDE for the Green's function of the two-dimensional wave equation,
having a product of Dirac deltas as its source.  The Green's function $G(x,y)$ {\em{per se}} is
$G(x,y)=iH_{0}^{(2)}(k\rho)/4,$
and it satisfies\footnote{We hasten to repeat that there is no claim whatsoever to originality at this point, the material
being in its entirety a part of common knowledge among physicists.}
\begin{eqnarray}
\left(\rule{0mm}{4mm}\nabla^{2}+k^{2}\right)G(x,y) & = & \delta(x)\delta(y)
\end{eqnarray}  
or else, following Fourier transformation
\begin{eqnarray}
\tilde{G}(w,y) & = & \int_{-\infty}^{\,\infty}
             e^{iwx}G(x,y)dx\,,
\end{eqnarray}
the reduced equation     
\begin{eqnarray}
\left[\frac{\partial^{2}}{\partial y^{2}}+k^{2}-w^{2}\right]\tilde{G}(w,y) & = & \delta(y)\,,
\end{eqnarray}
the solution of which latter is easily found to read
\begin{eqnarray}
\tilde{G}(w,y) & = & \frac{i}{2\sqrt{k^{2}-w^{2}\,}}\,e^{-i\sqrt{k^{2}-w^{2}\,}|y|}\,.
\end{eqnarray}
\newpage
\mbox{   }
\newline
\newline        
\newline
Then
\begin{eqnarray}
\int_{0}^{\infty}H_{0}^{(2)}(kx)\cos(wx)dx & = & -2i\tilde{G}(w,0)  \nonumber  \\
                                           & = & \frac{1}{\sqrt{k^{2}-w^{2}\,}}
\end{eqnarray}
and so we are done once again.  At the same time this latest intrusion by Dirac deltas compels us to acknowledge
that all three derivations of (1), neither Basset's nor ours in any way excluded, skate, in one mode or another,
upon a thin ice sheet of analytic delicacy.
\section{Remarks}
   It goes without saying that use of the alternate time dependence, $\exp(-i\omega t),\;\omega > 0,$ induces
complex conjugation {\em{en masse}},
$H_{0}^{(2)}(k\rho)\rightarrow H_{0}^{(1)}(k\rho),$ $\Im k \geq 0,$ $|\Im w|<\Im k$ as an initial strip of
analyticity, extendable thereafter to the whole plane when properly cut.

   On a somewhat different tack, one must remember that the limit $\beta\downarrow 0+$ cannot be recklessly enforced,
beginning with evaluation (10).  For example, joining (10) is the further elementary integral
\begin{eqnarray}
\int_{0}^{\,\infty}e^{-\beta x}\sin(x)dx & = & \frac{1}{1+\beta^{2}}
\end{eqnarray}
already implied in (15).  But, in connection with that limit, there is of course no intent
to set
\begin{eqnarray}
\int_{0}^{\,\infty}\sin(x)dx & = & 1  \nonumber  \\
\int_{0}^{\,\infty}\cos(x)dx & = & 0 \,,
\end{eqnarray}
both proposed quadratures being clearly nonsensical.  Meaningless likewise would be the tentative assignments
\begin{eqnarray}
\int_{0}^{\,\infty}\sin(x)\log(x)dx & = & C  \nonumber  \\
\int_{0}^{\,\infty}\cos(x)\log(x)dx & = & -\pi/2\,.
\end{eqnarray}

  A nagging anxiety that retaining in (9) only the lowest-order close-in Hankel terms may not be adequate
is alleviated somewhat by considering, and then at once dismissing yet another potential contribution proportional
to $\log(k/2w)$ [{\bf{16}}], {\em{viz.,}}
\begin{eqnarray}
\Delta & = & -\frac{2i}{\pi w}\log\left(\frac{k}{2w}\right)\sum_{r=1}^{\infty}(-1)^{r}
\left(\frac{k^{r}}{2^{r}w^{r}r!}\right)^{2} \lim{\atop_{\rule{-7.5mm}{2.5mm}\beta\,\downarrow\,0+}}
\int_{0}^{\,\infty}e^{-\beta x}x^{2r}\cos(x)dx\,.
\end{eqnarray}
\newpage
\mbox{   }
\newline
\newline        
\newline
However, with a view to (10),
\begin{eqnarray}
 \lim{\atop_{\rule{-7.5mm}{2.5mm}\beta\,\downarrow\,0+}}\int_{0}^{\,\infty}e^{-\beta x}x^{2r}\cos(x)dx & = &
 \lim{\atop_{\rule{-7.5mm}{2.5mm}\beta\,\downarrow\,0+}}\left(-\frac{d}{d\beta}\right)^{2r}        
 \sum_{n=0}^{\infty}(-1)^{n}\beta^{2n+1} \nonumber  \\
     & = & 0\,,
\end{eqnarray}    
and therefore the asymptotic development of $L(k,w)$ is entirely free of logarithmic dependence
at every power of multiplier $1/w.$  All other analytic fragments contributing to $H_{0}^{(2)}(kx/w)$
in the first line of (9) can be treated in similar fashion, with similarly null outcomes.

    And as for the constants $C$ and $D$ from (17), we can fix those via outright numerical quadrature of
\begin{eqnarray}
N(\eta) & = & \int_{0}^{\,\infty} e^{-\eta x} \sin(x)\log(x)dx
\end{eqnarray}
at, say, $\eta=\mu$ and $\eta=\nu,$ $\mu\neq\nu.$  Thus
\begin{equation}
\left(
\begin{array}{cc}
     1 & \mu \\
       &     \\
     1 & \nu
\end{array} \right)
\left(
\begin{array}{c}
     C  \\
        \\
     D
\end{array} \right)  = \left( \begin{array}{c}
                                (1+\mu^{2})N(\mu)+\frac{1}{2}\log(1+\mu^{2})+\mu\arctan(\mu) \\
                                                                                     \\
                                (1+\nu^{2})N(\nu)+\frac{1}{2}\log(1+\nu^{2})+\nu\arctan(\nu)
                              \end{array} \right)\,,                                                                                       
\end{equation}
the matrix on the left having $\nu-\mu\neq 0$ as determinant and hence being open to inversion as
\begin{equation}
\left(
\begin{array}{cc}
1 & \mu \\
  &     \\
1 & \nu
\end{array} \right)^{-1}=\frac{1}{\nu-\mu}\left(
\begin{array}{cc}
\nu & -\mu \\
    &      \\
 -1 &  1
\end{array} \right)\,.
\end{equation}
While system (31) can thus be solved in the usual way, one naturally imposes the self-consistency demand that
any other pair $\mu'\neq\nu'$ yield an identical outcome, apart from any inaccuracy attributable to the
numerical treatment of (30).  But, even though it be a {\em{sine qua non desideratum}}, we are in no
position to offer any sort of {\em{ab initio}} proof of its fulfillment.                  
%
\theendnotes
%
\section{References}
\parindent=0in

1.	Sommerfeld, A. (1896). Mathematische theorie der diffraction, {\em{Mathematische Annalen}}, 47:317-374.

2.	Sommerfeld, A. (2004). {\em{Mathematical Theory of Diffraction}} (English translation by Nagem, R. J.,
Zampolli, M., Sandri, G.), Springer Science+Business Media, New York, USA.

3.	Macdonald, H. M. (1902). {\em{Electric Waves, (Appendix D. Diffraction)}}, pp. 186-196, Cambridge University
Press, Cambridge, UK.

4.	Lamb, H. (1906). On Sommerfeld's diffraction problem; and on reflection by a parabolic mirror, {\em{Proc.
London Math. Soc.}}, 4: 190-203.

5.	Copson, E. T. (1946). On an integral equation arising in the theory of diffraction, {\em{Quart. J. Math.}}, 
(Oxford), 17: 19-34.

6.	Baker, B. B., Copson, E. T. (1949). {\em{The Mathematical Theory of Huygens' Principle (Second Edition)}},
Chapter V.  Diffraction by a Plane Screen, pp. 153-177, Oxford University Press, UK.

7.	Schwinger, J., DeRaad, Jr., L. L., Milton, K. M., Tsai, W-Y. (1998). {\em{Classical Electrodynamics}}, Chapter 48,
Diffraction II, pp. 509-521, Westview Press, Boulder, Colorado, USA.

8.	Jones, D. S. (1952). A simplifying technique in the solution of a class of diffraction problems, {\em{Quart. J.
Math.}}, 3(2): 189-196.

9.	Noble, B. (1958). {\em{Methods Based on the Wiener-Hopf Technique}}, Pergamon Press, London, UK.
\newpage
\mbox{   }
\newline

10.	Mittra, R., Lee, S. W. (1971). {\em{Analytical Techniques in the Theory of Guided Waves}}, The Macmillan
Company, New York, USA.

11. Clemmow, P. C. (1996). {\em{The Plane Wave Spectrum Representation of Electromagnetic Fields}}, IEEE
Press, Piscatawaty, New Jersey, USA.

12. Born, M., Wolf, (1999). {\em{Principles of Optics (Seventh \{Expanded\} Edition)}}, Chapter XI. Rigorous
diffraction theory (contributed by Clemmow, P. C.), pp. 667-673, Cambridge University Press, Cambridge, UK.

13. Carrier, G. F., Krook, M., Pearson, C. E. (2005). {\em{Functions of a Complex Variable (Theory and
Technique)}}, Chapter 8. Special Techniques, pp. 376-408, Society for Industrial and Applied Mathematics,
Philadelphia, USA.

14.	Basset, A. B. (1988). {\em{A Treatise on Hydrodynamics, Volume II}}, p. 32, Deighton, Bell and Co., Cambridge, UK.
 
15.	Watson, G. N. (1944). {\em{A Treatise on the Theory of Bessel Functions, (Second Edition)}}, p. 388, Cambridge, UK.

16.	Copson, E. T. (1935). {\em{An Introduction to the Theory of Functions of a Complex Variable}}, p. 329, Oxford, UK.


\end{document}